\documentclass[amsfonts,twoside]{article}
\usepackage{amsmath,amsfonts,amscd,amsthm,amssymb,enumerate,graphics,graphicx,tabularx,array}

 \font\chuto=cmbx10 at 14pt

\font\ita=cmti9

\font\ita=cmti9   

\newtheorem{thm}{Theorem}[subsection]
\newtheorem{lem}[thm]{Lemma}
\newtheorem{prop}[thm]{Proposition}

\theoremstyle{definition}
\newtheorem{defn}[thm]{Definition}
\newtheorem{ex}[thm]{Example}
\newtheorem{rem}[thm]{Remark}

\title{
{\chuto The Second Cohomology Group of Elementary Quadratic Lie Superalgebras and Classifying a Subclass of 8-dimensional Solvable Quadratic Lie Superalgebras}}

\author{{\bf Cao Tran Tu Hai$^*$, Duong Minh Thanh$^\dagger$ , 
Le Anh Vu$^{**}$}}

\begin{document}
\date{}
\maketitle
\label{firstpage}

\begin{center}
	\footnotesize{\ita  $^*$  Le Quy Don High School, Ninh Thuan, Vietnam}\\
	\footnotesize{\it Email: tuhai.thptlequydon@ninhthuan.edu.vn}
\vskip0.2cm	
	\footnotesize{\ita  ${}^\dagger$Ho Chi Minh city University of Education, Vietnam}\\
	\footnotesize{\it Email: thanhdmi@hcmup.edu.vn}
\vskip0.2cm	
	\footnotesize{\ita  ${}^{**}$University of Economics and Law, VNU-HCMC, Vietnam}\\
	\footnotesize{\it Email: vula@uel.edu.vn}
\end{center}

\vskip0.8cm

\begin{abstract}
By definition, a quadratic Lie superalgebra is a Lie superalgebra endowed with a non-degenerate supersymmetric bilinear form which satisfies the even and invariant properties. In this paper we calculate all of the second cohomology group of elementary quadratic Lie superalgebras which have been classified in \cite{DU14} by applying the super-Poisson bracket on the super exterior algebra. Besides, we give the classification of 8-dimensional solvable quadratic Lie superalgebras having 6-dimensional indecomposable even part. The method is based on the double extension and classification results of adjoint orbits of the Lie algebra $\mathfrak{s}\mathfrak{p}(2)$.
\end{abstract}

\vskip0.8cm




\section*{Introduction}
As far as we know, the Killing form of a Lie superalgebra is supersymmetric, invariant and even. In some special cases, it also satisfies the non-degeneracy. Those  lead to study of Lie superalgebras endowed with a supersymmetric, invariant, even and non-degenerate bilinear form. Such Lie superalgebras are called quadratic Lie superalgebras.

Consider the constructive aspect, a non-trivial quadratic Lie algebra could be considered as a double extension (a combination of central extension and semi-direct product) of a quadratic Lie algebra of smaller dimension (see \cite{MR85}). Moreover, every solvable quadratic Lie algebra is also isometrically isomorphic to either a T*-extension of a certain Lie algebra with its dual space or an ideal of codimension one of a T*-extension (see \cite{Bor97}). Both of these conceptions were generalized by S. Bajo, H. Benamor, S. Benayadi and M. Bordemann for quadratic Lie superalgebras in \cite{BB97} and \cite{BBB}. In what follows, it is of interest to research algebras endowed with an invariant and non-degenerate bilinear form as well as their applications.

Another concerned problem is to describe the cohomology of Lie superalgebras, which is an important tool in mathematrics and theoretical physics. 
A classical example of a constant such that 
D. B. Fuchs and D. A. Leites in \cite{FL84} calculated the cohomology groups of the classical Lie superalgebras with trivial coefficients,
Y. C. Su and R. B. Zhang in \cite{SZ07} computed explicitly the first and second cohomology groups of the classical Lie superalgebras $\mathfrak{sl}_{m|n}$ and $\mathfrak{osp}_{2|2n}$ with coefficients in the finite-dimensional irreducible modules and the Kac modules, W. Bai and W. Liu in \cite{BL17} described the cohomology groups of Heisenberg Lie superalgebras.

We shall now describe a situation in which
$\mathfrak{g}$ is a quadratic Lie superalgebra, then the algorithm of describing the second cohomology groups with coefficients in  $\mathbb{C}$ of $\mathfrak{g}$ is relevant to describing the super-Poisson bracket of the super-exterior algebra of $\mathfrak{g}$.
In a consequence of this result,
 we can list all of their one-dimensional double extensions. This provides much information for the classification of quadratic Lie superalgebras. Our goal in this article is to calculate the second cohomology group of all elementary quadratic Lie superalgebras classified in \cite{DU14} and give a classification of 8-dimensional solvable quadratic Lie superalgebras having 6-dimensional indecomposable even part. The classification is based on the double extension and the classification of adjoint orbits of the Lie algebra $\mathfrak{s}\mathfrak{p}(2)$.

The paper will be organized as follows: The first section is devoted to recall some basic concepts and results of Lie superalgebras, cohomology of Lie superalgebras and quadratic Lie superalgebras. The second section gives the second cohomology group of all elementary quadratic Lie superalgebras classified in \cite{DU14} by using the super $\mathbb{Z}\times {\mathbb{Z}_2}-$ Poisson bracket in the super-exterior algebra. The last section gives the classification of 8-dimensional solvable quadratic Lie superalgebras having 6-dimensional indecomposable even part by applying the results of a double extension of quadratic Lie superalgebras in \cite{BB97}, \cite{BBB} and the classification of adjoint orbits of the Lie algebra $\mathfrak{sp}(2)$. 

All vector spaces considered in throughout the paper are finite-dimensional complex vector spaces.

\section{Cohomology of Lie Superalgebras and \\Quadratic Lie Superalgebras}

In this section, we recall some preliminary concepts and basic results which will be used later. For details we refer the reader to the paper \cite{FL84} of D.B.Fuchs, D.A.Leites and the paper \cite{PU07} of G. Pinczon, R. Ushirobira.

\subsection{Lie Superalgebras and Cohomology}

\begin{defn}
{\it A Lie superalgebra} $\mathfrak{g}$ is a $\mathbb{Z}_{2}-$graded vector space $ \mathfrak{g} = \mathfrak{g}_{\overline 0}  \oplus \mathfrak{g}_{\overline 1}$ endowed with a Lie super bracket [.,.] that satisfies the following conditions: 
\begin{enumerate} [(i)] 
\item	
		The Lie super bracket [.,.] is bilinear and $[ \mathfrak{g}_{x}, \mathfrak{g}_{y} ] \subset \mathfrak{g}_{x+y}$\, ({\it grading});
\item 
 		$ [X,Y] = -(-1)^{xy} [Y,X]$ \,({\it skew-supersymmetry});
\item 
		 $(-1)^{zx} \left[  [X,Y], Z \right]+ (-1)^{xy} \left [ [Y,Z],X  \right]+ (-1)^{yz} \left[ [Z,X],Y \right]=0$ ({\it  super Jacobi identity})
\end{enumerate}
for all $x, y, z \in \mathbb{Z}_{2}$, $X\in \mathfrak{g}_{x},Y\in \mathfrak{g}_{y}, Z\in  \mathfrak{g}_{z}$.
\end{defn}

 \begin{defn}
Let 
$\mathfrak{g}=\mathfrak{g}_{\overline 0} \oplus \mathfrak{g}_{\overline 1}$ 
be a Lie superalgebra. Denote by 
$Alt(\mathfrak{g}_{\overline 0}, \mathbb{C})$
 the algebra of alternating multilinear forms on 
 $\mathfrak{g}_{\overline 0}$ 
  and by 
 $Sym(\mathfrak{g}_{\overline 1}, \mathbb{C})$ 
  the algebra of symmetric multilinear forms on 
 $\mathfrak{g}_{\overline 1}$. 
 We define a 
 $\mathbb{Z}\times \mathbb{Z}_{2}$-gradation on 
 $Alt(\mathfrak{g}_{\overline 0},\mathbb{C})$ 
 and on 
 $Sym(\mathfrak{g}_{\overline 1},\mathbb{C})$ by
$$Alt^{(i,\overline 0)}(\mathfrak{g}_{\overline 0},\mathbb{C})= Alt^{i}(\mathfrak{g}_{\overline 0},\mathbb{C}),Alt^{(i,\overline 1)}(\mathfrak{g}_{\overline 0},\mathbb{C})=
\{ 0 \}$$
and
$$Sym^{(i,\overline i)}(\mathfrak{g}_{\overline 1},\mathbb{C})=Sym^{i}(\mathfrak{g}_{\overline 1},\mathbb{C}),Sym^{(i,\overline{j})}(\mathfrak{g}_{\overline 1},\mathbb{C})=
\{0\}$$
where $i,j\in \mathbb{Z}$; $\overline{i},\overline{j}\in \mathbb{Z}_2$
 are respectively the residue classes modulo 2 of $i, j$ and $\overline{i}\ne \overline{j}$.
The super-exterior algebra of $\mathfrak{g}$ is 
$C(\mathfrak{g},\mathbb{C})=Alt(\mathfrak{g}_{\overline 0},\mathbb{C})\otimes Sym(\mathfrak{g}_{\overline 1},\mathbb{C})$ 
endowed with the super-exterior product on $C(\mathfrak{g},\mathbb{C})$ defined by
$$\left( \Omega \otimes F \right)\wedge \left( \Omega '\otimes F' \right)= (-1)^{f\omega '}\left( \Omega \wedge \Omega ' \right)\otimes FF',$$
for all $\Omega \in Alt(\mathfrak{g}_{\overline 0},\mathbb{C}),\Omega '\in Alt^{\omega '}(\mathfrak{g}_{\overline 0},\mathbb{C}),F\in Sym^{f}(\mathfrak{g}_{\overline 1},\mathbb{C}),F'\in Sym(\mathfrak{g}_{\overline 1},\mathbb{C})$.
\end{defn}
Remark that  $C(\mathfrak{g},\mathbb{C})$ is a $\mathbb{Z}\times \mathbb{Z}_2-$graded algebra. More precisely, in terms of $\mathbb{Z}-$gradation, one has
$$C^n(\mathfrak{g},\mathbb{C})=\underset{m=0}{\overset{n}{\mathop{\oplus }}}\,\left( Alt^m (\mathfrak{g}_{\overline 0},\mathbb{C})\otimes Sym^{n-m} (\mathfrak{g}_{\overline 1},\,\mathbb{C}) \right),C^0 (\mathfrak{g},\mathbb{C})=\mathbb{C},$$ 
and in terms of  $\mathbb{Z}_{2}-$gradation, 
$$C_{\overline 0}(\mathfrak{g},\mathbb{C})=Alt(\mathfrak{g}_{\overline 0},\mathbb{C})\otimes \left( \underset{j\ge 0}{\mathop{\oplus }}\,Sym^{2j}(\mathfrak{g}_{\overline 1},\mathbb{C}) \right)$$
$$ \text{ and } C_{\overline 1}(\mathfrak{g},\mathbb{C})=Alt(\mathfrak{g}_{\overline 0},\mathbb{C})\otimes \left( \underset{j\ge 0}{\mathop{\oplus }}\,Symm^{2j+1}(\mathfrak{g}_{\overline 1},\mathbb{C}) \right).$$

\begin{defn}
Denote by $End(C(\mathfrak{g},\mathbb{C}))$ the space of endomorphisms on $C(\mathfrak{g},\mathbb{C})$. A homogeneous endomorphism $D\in End(C(\mathfrak{g},\mathbb{C}))$ of degree $(n,d)$ is called {\it a superderivation} of  $C(\mathfrak{g},\mathbb{C})$ if it satisfies the following condition:
$$D(A\wedge A')=D(A)\wedge A'+{(-1)^{na+db}}A\wedge D(A), \forall A\in C^{(a,b)}(\mathfrak{g},\mathbb{C}),\forall A'\in C(\mathfrak{g},\mathbb{C}).$$
\end{defn}
Denote by $Der_{d}^{n}(C(\mathfrak{g},\mathbb{C}))$ the space of superderivations of degree $(n,d)$ of  $C(\mathfrak{g},\mathbb{C})$. Then we have a $\mathbb{Z}\times {\mathbb{Z}_2}-$gradation of the space of superderivations of $C(\mathfrak{g},\mathbb{C})$:
$$Der(C(\mathfrak{g},\mathbb{C}))=\underset{(n,d)\in \mathbb{Z}\times {\mathbb{Z}_{2}}}{\mathop{\oplus }}\, Der_{d}^{n}(C(\mathfrak{g},\mathbb{C})).$$

\begin{ex}
Let $X\in {\mathfrak{g}_{x}}$ be a homogeneous element in $\mathfrak{g}$ of degree $x$ and define the endomorphism $i_X$ of   
$C(\mathfrak{g},\mathbb{C})$ by  
$$i_X (A)\left( X_1,\ldots ,X_{a-1} \right)={(-1)}^{xb} A\left( X,X_1,\ldots ,X_{a-1} \right)$$
 for all $A\in C^{(a,b)}(\mathfrak{g},\mathbb{C}); X_1,\ldots, X_{a-1} \in \mathfrak{g}.$ 
Then 
$$i_X(A\wedge A')=i_X (A)\wedge A'+(-1)^{-a+xb}  A\wedge i_X (A)$$
 for all $A\in C^{(a,b)}(\mathfrak{g},\mathbb{C})$, $A'\in C(\mathfrak{g},\mathbb{C}).$ It means that $i_X$ is a superderivation of  degree $(-1,x)$. 
\end{ex}

Given $k\ge 0$, the differential operator $\delta _{k}:C^{k}(\mathfrak{g},\mathbb{C})\to C^{k+1}(\mathfrak{g},\mathbb{C})$ is a superderivation of degree $\left( 1,\overline{0} \right)$ defined by
 \begin{multline*}
\delta _{k} \omega \left( X_0 ,\ldots , X_{k} \right)\\
=\sum\limits_{r<s} {(-1)}^{{s} + x_{s} \left( x_{r+1} + \ldots + x_{s-1} \right)} \omega \left( X_{0}, \ldots , X_{r-1}, \left[ X_{r}, X_{s} \right], X_{r+1}, \ldots ,\widehat{X_{s}},\ldots ,X_{k} \right),
\end{multline*}
for all $\omega \in C^{k}(\mathfrak{g},\mathbb{C})$, $ X_{0} \in \mathfrak{g}_{{x}_{0}}, \ldots , X_{k} \in \mathfrak{g}_{x_{k}}$, where the sign $\widehat{X_{s}}$ indicates that  the element $X_{s}$ is omitted. It is easy to check that $\delta ^{2} =\delta \circ \delta =0.$ By convention, $\delta _{0} =0.$

An element $\omega \in C^{k}(\mathfrak{g},\mathbb{C})$ is called a $k$-{\it cocycle} if  $\delta _{k}\omega =0$ or a $k$-{\it coboundary } if  there exists $\varphi \in C^{k-1}(\mathfrak{g},\mathbb{C})$ such that $\omega = \delta _{k-1}\varphi $.

We denote by $Z^{k} (\mathfrak{g},\mathbb{C})$ the set of all $k$-cocycles and by $B^{k}(\mathfrak{g},\mathbb{C})$ the set of all $k$-coboundaries. That is $Z^{k}(\mathfrak{g},\mathbb{C})= Ker \delta _{k}$ and $B^{k}(\mathfrak{g},\mathbb{C})=Im \delta _{k-1}$. Clearly, $B^{k} (\mathfrak{g},\mathbb{C})\subset Z^{k}(\mathfrak{g},\mathbb{C})$.  The quotient space  $Z^{k}(\mathfrak{g},\mathbb{C})/B^{k}(\mathfrak{g},\mathbb{C})$ is denoted by $H^{k}(\mathfrak{g},\mathbb{C})$ and called the $k$ {\it -cohomology groups} of $\mathfrak{g}$ with trivial coefficients. 

\begin{defn}
 The dimension of the $k$-cohomology group $H^{k}(\mathfrak{g},\mathbb{C})$ is called  the {\it $k$-th Betti number} of $\mathfrak{g}$ and denoted by $b_{k}(\mathfrak{g})$.
\end{defn}

\begin{ex}
 (see \cite{BL17}) Let the  Heisenberg Lie superalgebra  
 $$\mathfrak{h}_{2n+1,m}
 =\mathfrak{g}_{\overline{0}} \oplus \mathfrak{g}_{\overline{1}}
 =\mathbb{C} \{  Z,X_{1}, \ldots , X_{n}, X_{n+1}, \ldots , X_{2n} \} \oplus \mathbb{C} \{  Y_{1}, \ldots ,Y_{m} \}$$
with non-zero super brackets 
$$ [ X_{i}, X_{n+i} ] =Z, [ Y_{j}, Y_{j}] = Z,\forall i=\overline{1,n},j=\overline{1,m}.$$ 
It is easy to compute that $\delta X_{i}^{*}=\delta Y_{j}^{*}=0,$ for all $i=\overline{1,2n},j=\overline{1,m}$ and 
$$\delta Z^{*}=\sum\limits_{i=1}^{n}{X_{n+i}^{*}\wedge X_{i}^{*}}-\frac{1}{2}\sum\limits_{j=1}^{m}{Y_{j}^{*}Y_{j}^{*}}.$$
For the second cohomology group, we have $\delta \left( Z^{*}\wedge \omega  \right)=\delta Z^{*}\wedge \omega -Z^{*}\wedge \delta \omega =0$ if and only if $\omega =0$. Then
$$Z^{2}(\mathfrak{h}_{2n+1,m},\mathbb{C})=\left\{ X_{i}^{*}\wedge X_{j}^{*},X_{i}^{*}\otimes Y_{k}^{*},Y_{k}^{*} Y_{l}^{*}:i,j=\overline{1,2n},i\ne j,k,l=\overline{1,m} \right\}$$
$${\text { and }} \dim{Z^{2}}(\mathfrak{h}_{2n+1,m},\mathbb{C})=\left( 
\begin{matrix}
   2n  \\
   2  \\
\end{matrix} \right)+2n.m+m+\left( \begin{matrix}
   m  \\
   2  \\
\end{matrix} \right)=2n^{2} -n + 2nm + \frac{m^{2}+m}{2}.$$
Accordingly, $b_{2}(\mathfrak{h}_{2n+1,m})=2n^{2} - n + 2nm + \dfrac{m^{2}+m}{2} -1.$
\end{ex}

\subsection{Quadratic Lie Superalgebras}

\begin{defn}
 Let $\mathfrak{g}=\mathfrak{g}_{\overline 0} \oplus \mathfrak{g}_{\overline 1}$ be a Lie superalgebra.  Assume that $B$ is a bilinear form defined on $\mathfrak{g}$ such that it satisfies the following properties:
\begin{enumerate} 
\item [(i)]	$B(X,Y)=(-1)^{xy} B(Y,X)$, $\forall X\in \mathfrak{g}_{x},Y\in \mathfrak{g}_{y}$ ({\it supersymmetric});
\item [(ii)]	$B\left( \left[ X,Y \right],Z \right)=B\left( X,\left[ Y,Z \right] \right)$ for all $X,Y,Z\in \mathfrak{g}$ ({\it invariant});
\item [(iii)]	$B(X,Y)=0$, $\forall Y\in \mathfrak{g}$ implies $X=0$ ({\it non-degenerate}).
\end{enumerate}
The pair ($\mathfrak{g}, B$) is called {\it a quadratic Lie superalgebra} if $B$ is even, that is  
$$B(X,Y) = 0; \, \forall X \in \mathfrak{g}_{\overline 0}, Y \in \mathfrak{g}_{\overline 1}.$$ 
In this case, it is easy to check that $\left(\mathfrak{g}_{\overline 0},B\left| _{\mathfrak{g}_{\overline 0}}\times \mathfrak{g}_{\overline 0} \right. \right)$ is a quadratic Lie algebra and $\left( \mathfrak{g}_{\overline 1},B\left| _{\mathfrak{g}_{\overline 1}}\times \mathfrak{g}_{\overline 1} \right. \right)$ is a $\mathfrak{g}_{\overline 0}-$module endowed with a symplectic structure.
\end{defn}

Let $\left( \mathfrak{g},B \right), \left( \mathfrak{g}',B' \right)$ be two quadratic Lie superalgebras. We say $\left( \mathfrak{g},B \right)$ and $\left( \mathfrak{g}',B' \right)$ {\it isometrically isomorphic} (or i-isomorphic, for short) if there exists a Lie superalgebra isomorphism $A$ from $\mathfrak{g}$ onto $\mathfrak{g}'$ satisfying 
$$B'\left( A(X),A(Y) \right)=B\left( X,Y \right), \forall X,Y\in \mathfrak{g}.$$
 Then $A$ is called an {\it i-isomorphism}. We write  $\left( \mathfrak{g},B \right)\overset{i}{\mathop{\cong }}\,\left( \mathfrak{g}',B' \right)$.

\begin{defn} Let $\left( \mathfrak{g},B \right)$ be a quadratic Lie superalgebra and $\Im $ be a graded ideal of $\mathfrak{g}$. 
\begin{enumerate} [(i)]
\item 
		$\Im $ is called {\it non-degenerate} if  the restriction of $B$ to $\Im \times \Im $ is non-degenerate. Otherwise, we say $\Im $ {\it degenerate}.
\item 
		$\left( \mathfrak{g},B \right)$ is called {\it irreducible} if $\mathfrak{g}$ does not have any non-degenerate graded ideal excepting $ \{ 0 \}$ and $\Im $.
\item
		 A non-degenerate ideal $\Im $ is called irreducible if $\Im $ does not have any non-degenerate graded ideal excepting $ \{ 0 \}$ and $\Im $.
\item
		 Ideal $\Im $ is called {\it totally isotropic} if $B\left(\Im ,\Im \right)=\{0 \}$.
\end{enumerate}
\end{defn}

The following proposition reduces the study of quadratic Lie superalgebras to non-degenerate graded ideals.

\begin{prop}[see  \cite{BBB}]

Let $\left( \mathfrak{g},B \right)$ be a quadratic Lie superalgebra and $\Im $ be a graded ideal of $\mathfrak{g}$. Then  $\Im^{\bot}$ is also a graded ideal of $\mathfrak{g}$. In addition, if  $\Im $ is non-degenerate then so is $\Im ^{\bot }$, $\left[ \Im ,\Im ^{\bot} \right]=\left\{ 0 \right\}$ and $\Im \cap \Im ^{\bot } =\left\{ 0 \right\}$. In this case, we denote $\mathfrak{g} = \Im \overset{\bot}{\mathop{\oplus }}\, \Im ^{\bot }$. \hfill $\square$
\end{prop}

\section{The Second Cohomology Group of \\Elementary Quadratic Lie Superalgebras}

In this section, we will compute the second cohomology group of all elementary quadratic Lie superalgebras classified in \cite{DU14}. Firstly, we recall the concept of the super $\mathbb{Z} \times \mathbb{Z}_{2}-$ Poisson bracket on the super-exterior algebra of a quadratic Lie superalgebra, which is used to give a new way of description of cohomology. 

\subsection{The Super $\mathbb{Z}\times \mathbb{Z}_{2}-$Poisson Bracket on The Super-exterior Algebra} 
Let  $\mathfrak{g}= \mathfrak{g}_{\overline 0}\oplus \mathfrak{g}_{\overline 1}$ be a  $\mathbb{Z}_{2}-$graded vector space equipped with a non-degenerate even supersymmetric bilinear  form $B$. In this case, $\mathfrak{g}_{\overline 1}$ is a symplectic vector space. Hence, the dimension of $\mathfrak{g}_{\overline{1}}$ must be even and  $\mathfrak{g}$ is aslo called a {\it quadratic} $\mathbb{Z}_{2}-$graded vector space. Now we recall the definition of the Poisson bracket on $Sym(\mathfrak{g}_{\overline 1})$ and the super-Poisson bracket on
 $Alt(\mathfrak{g}_{\overline 0},\mathbb{C})$ which are used later. 
 
 Let 
  $\left\{ X_{1}, \ldots , X_{n},Y_{1}, \ldots , Y_{n} \right\}$ 
  be a Darboux basis of  $\mathfrak{g}_{\overline 1}$, i.e. we have 
  $$B(X_{i}, X_{j})=B(Y_{i}, Y_{j})=0, B(X_{i}, Y_{j}) = \delta _{ij},$$ 
  for all 
 $i,j=\overline{1,n}$. Let $\{ p_{1}, \ldots , p_{n},q_{1}, \ldots , q_{n} \}$ 
 be its dual basis. Then the algebra 
 $Sym(\mathfrak{g}_{\overline 1}, \mathbb{C})$ 
 regarded as the polynomial algebra 
 $\mathbb{C} [p_{1}, \ldots ,p_{n},q_{1}, \ldots ,q_{n}]$ 
 is equipped with the Poisson bracket as follows:
$$\left\{ F,G \right\}= \sum\limits_{i=1}^{n} {\left( \frac{\partial F}{\partial p_{i}}\frac{\partial G}{\partial q_{i}}-\frac{\partial F}{\partial q_{i}}\frac{\partial G}{\partial p_{i}} \right)},\forall F,G\in Sym({\mathfrak{g}}_{\overline{1}},\mathbb{C}).$$
For any $X\in \mathfrak{g}_{\overline 0}$, let $i_{X}$ be the derivation of 
$Alt(\mathfrak{g}_{\overline 0},\mathbb{C})$ 
defined by 
 $$\iota _{X} (\Omega )(Z_{1}, \ldots , Z_{k})=\Omega (X,Z_{1}, \ldots , Z_{k}), 
\forall \Omega \in Alt^{k+1} (\mathfrak{g}_{\overline 0}, \mathbb{C}),X, Z_{1}, \ldots , Z_{k} \in \mathfrak{g}_{\overline 0}, k\ge 0$$ 
and $\iota _{X}(1)=0$. Let 
$\{ Z_{1}, \ldots , Z_{m} \}$ 
be a fixed orthonormal basis of 
$\mathfrak{g}_{\overline 0}$. 
The super-Poisson bracket on 
$Alt(\mathfrak{g}_{\overline 0},\mathbb{C})$ 
is defined  by  (see \cite{PU07}):
$$\left\{ \Omega ,\Omega ' \right\}={(-1)^{k+1}}\sum\limits_{j=1}^{m}{ \iota _{Z_{j}} (\Omega )\wedge \iota _{Z_j} (\Omega ')},\forall \Omega \in Alt^{k}(\mathfrak{g}_{\overline 0}, \mathbb{C}), \Omega ' \in Alt(\mathfrak{g}_{\overline 0},\mathbb{C}).$$
Next, the super $\mathbb{Z}\times \mathbb{Z}_{2}-$Poisson bracket on $C(\mathfrak{g},\mathbb{C})$ is given by:
 $$\left\{ \Omega \otimes F,\Omega '\otimes G \right\}= (-1)^{f \omega '} \left( \left\{ \Omega ,\Omega ' \right\}\otimes FG + \Omega \wedge \Omega ' \otimes \left\{ F,G \right\} \right);$$
for any $\Omega \in Alt(\mathfrak{g}_{\overline 0}, \mathbb{C}),\Omega '\in Alt^{\omega '}(\mathfrak{g}_{\overline 0},\mathbb{C}), F \in Sym^{f} ( \mathfrak{g}_{\overline 1}, \mathbb{C}),G \in Sym( \mathfrak{g}_{\overline 1}, \mathbb{C}).$

\begin{prop}[see \cite{DU14}, \cite{PU07}]
The algebra $C(\mathfrak{g},\mathbb{C})$ is a graded Lie algebra with the super $\mathbb{Z} \times \mathbb{Z}_{2}-$Poisson bracket. More precisely, for all $A \in C^{(a,b)} ( \mathfrak{g},\mathbb{C}), A' \in C^{(a',b')} (\mathfrak{g},\mathbb{C})$ and $A'' \in {C}(\mathfrak{g},\mathbb{C})$, one has
\begin{enumerate}[(i)]
\item 
	$\{A',A\}= -{(-1)^{aa'+bb'}}\{A,A'\},$
\item 
	$\{ \{ A, A' \},A'' \} =\{A,\{ A',A'' \} \}-{(-1)^{aa'+bb'}}\{ A',\{ A,A'' \} \}.$
\end{enumerate} 
Furthermore, $\{ A,A'\wedge A''\} = \{ A,A' \} \wedge A''+{(-1)^{aa'+bb'}}A'\wedge \{ A,A'' \}.$ \hfill $\square$
\end{prop}

\noindent Now, we choose an arbitrary basis 
$\left\{ X_{\overline 0}^{1}, \ldots ,X_{\overline 0}^{m} \right\}$ of 
$\mathfrak{g}_{\overline 0}$.
Denote by $\{ \alpha _{1},\ldots , \alpha _{m} \}$ its dual basis. Let  
$\left\{ Y_{\overline 0}^{1}, \ldots,Y_{\overline 0}^{m} \right\}$ 
be the basis of  $ \mathfrak{g}_{\overline 0}$ defined by $\alpha _{i} = B(Y_{\overline 0}^{i},\centerdot ).$ That means 
$$B(Y_{\overline 0}^i, X_{\overline 0}^j) = \delta _{ij}; \, \forall i, j = 1, \ldots , n.$$
Set 
$\left\{ X_{\overline 1}^{1},\ldots ,X_{\overline 1}^{n}, Y_{\overline 1}^{1},\ldots ,Y_{\overline 1}^{n} \right\}$ 
be a Darboux basis of  
$\mathfrak{g}_{\overline 1}$. 
Then the super $\mathbb{Z}\times \mathbb{Z}_{2}-$Poisson bracket on 
$C(\mathfrak{g},\mathbb{C})$ is also given by 
\begin{multline*}
\left\{ A,A' \right\}=(-1)^{\omega +f+1} \sum\limits_{i,j=1}^{m}{B(Y_{\overline 0}^{i},Y_{\overline 0}^{j}).\iota _{X_{\overline 0}^{i}} (A) \wedge \iota _{X_{\overline 0}^{j}} (A')}\\
+(-1)^{\omega } \sum\limits_{k=1}^{n}{\left( \iota _{X_{\overline 1}^{k}} (A)\wedge \iota _{Y_{\overline 1}^{k}} (A')- \iota _{Y_{\overline 1}^{k}} (A)\wedge \iota _{X_{\overline 1}^{k}}(A') \right)},
\end{multline*}
for all $A\in Alt^{\omega } (\mathfrak{g}_{\overline 0}, \mathbb{C})\otimes Sym^{f} (\mathfrak{g}_{\overline 1},\mathbb{C}),A'\in C(\mathfrak{g},\mathbb{C})$ (see \cite{DU14}).

\begin{rem}
 In \cite{DU14}, the authors introduced a useful tool. That was the 3-form $I$ defined on any quadratic Lie superalgebra $(\mathfrak{g},B)$ as follows
$$I\left( X,Y,Z \right)=B([X,Y],Z),\forall X,Y,Z\in \mathfrak{g}.$$
This 3-form is called the {\it 3-form associated }to  $\mathfrak{g}$.  It is easy to prove that $I$ is the homogeneous element of degree 
$(3,\overline{0})$ 
in the $\mathbb{Z}\times \mathbb{Z} _{2}-$graded algebra 
$C(\mathfrak{g},\mathbb{C})=Alt( \mathfrak{g}_{\overline 0}, \mathbb{C})\otimes Sym(\mathfrak{g}_{\overline 1}, \mathbb{C})$, $\left\{ I,I \right\}=0$ and $\delta =-\left\{ I,. \right\}$ (see \cite{DU14}, Proposition 1.11). Using this proposition, the cohomology group  
$H^{k} (\mathfrak{g},\mathbb{C})$ can be computed through the super 
$\mathbb{Z}\times \mathbb{Z}_{2}-$Poisson bracket.
\end{rem}

\subsection{Elementary Quadratic Lie Superalgebras}\label{2.2}

The main result of Section 2 is the description of the second cohomology group of elementary quadratic Lie superalgebras which have classified in \cite{DU14}. There are exactly three superalgebras as follows.
\begin{itemize}
	\item[\bf {2.2.1.}] $\mathfrak{g}_{4,1}^{s}$ $=\left( \mathbb{C} X_{\overline 0} \oplus \mathbb{C} Y_{\overline 0} \right) \oplus \left( \mathbb{C} X_{\overline 1} \oplus \mathbb{C} Y_{\overline 1} \right)$,
	where 
        $\mathfrak{g}_{\overline 0} = span\{ X_{\overline 0}, Y_{\overline 0} \}, 
        \mathfrak{g}_{\overline 1} = span\{ X_{\overline 1}, Y_{\overline 1} \}$. 
The bilinear form $B$ is defined by 
 $B\left(X_{\overline 0}, Y_{\overline 0} \right)=1, 
 B\left( X_{\overline 1}, Y_{\overline 1} \right)=1$, 
the others are zero and the Lie super bracket is given by 
      $$\left[ Y_{\overline 1}, Y_{\overline 1} \right] = - 2X_{\overline 0}, 
      \left[ Y_{\overline 0}, Y_{\overline 1} \right] = - 2X_{\overline 1}.$$
	\item[\bf {2.2.2.}] $\mathfrak{g}_{4,2}^{s} =\left( \mathbb{C} X_{\overline 0} \oplus \mathbb{C} Y_{\overline 0} \right)\oplus \left( \mathbb{C}X_{\overline 1} \oplus \mathbb{C} Y_{\overline 1} \right)$,
	where 
$ \mathfrak{g}_{\overline 0} =span\{ X_{\overline 0}, Y_{\overline 0} \}, \mathfrak{g}_{\overline 1} = span\{ X_{\overline 1}, Y_{\overline 1} \}$. 
The bilinear form $B$ is defined by 
$B\left( X_{\overline 0}, Y_{\overline 0} \right)=1, B\left( X_{\overline 1}, Y_{\overline 1} \right)=1$, 
the others are zero and  the Lie super bracket is given by 
$$\left[ X_{\overline 1}, Y_{\overline 1} \right] = X_{\overline 0}, \left[ Y_{\overline 0}, X_{\overline 1} \right] = X_{\overline 1}, \left[ Y_{\overline 0}, Y_{\overline 1} \right] = - Y_{\overline 1}.$$
	\item[\bf {2.2.3.}] $\mathfrak{g}_{6}^{s} =\left( \mathbb{C} X_{\overline 0} \oplus \mathbb{C} Y_{\overline 0} \right) \oplus \left( \mathbb{C} X_{\overline 1} \oplus \mathbb{C} Y_{\overline 1} \oplus \mathbb{C} Z_{\overline 1} \oplus \mathbb{C} T_{\overline 1} \right)$,
	where 
$\mathfrak{g}_{\overline 0} = span\{ X_{\overline 0}, Y_{\overline 0} \}$, 
$\mathfrak{g}_{\overline 1} = span\{ X_{\overline 1}, Y_{\overline 1}, Z_{\overline 1},  T_{\overline 1} \}$. 
The bilinear form $B$ is defined by
 $B\left( X_{\overline 0}, Y_{\overline 0} \right)=1$,
  $B\left( X_{\overline 1}, Z_{\overline 1} \right)=1$, 
  $B\left( Y_{\overline 1}, T_{\overline 1} \right)=1$, 
  the others are zero and the Lie super bracket is given by  
       $$\left[ Z_{\overline 1}, T_{\overline 1} \right] = - X_{\overline 0},
       \left[ Y_{\overline 0}, Z_{\overline 1} \right] = - Y_{\overline 1}, 
       \left[ Y_{\overline 0}, T_{\overline 1} \right] = - X_{\overline 1}.$$
\end{itemize}

\subsection{The Second Cohomology Group of Elementary \\ Quadratic Lie Superalgebras} \label{2.3}

Now we will introduce the first result of the paper. Namely, we will describe the second cohomology group of the elementary quadratic Lie superalgebras which have listed in Subsection 2.2. 

\begin{thm} \label{secondcohomology}
	With notations being as above in Subsection \ref{2.2}, the second cohomology group of the elementary quadratic Lie superalgebras are described as follows
	\begin{enumerate}	 [(i)]
		\item 
		$H^2(\mathfrak{g}_{4,1}^s,\mathbb{C}) 
		= span\left\{ 
      		          \left[ Y_{\overline 0}^{*}\otimes X_{\overline 1}^{*} \right],
		               \left[ X_{\overline 1}^{*}Y_{\overline 1}^{*}
		                            -2X_{\overline 0}^{*}\wedge Y_{\overline 0}^{*}\right]
		       \right\}$ 
		      where  
		      $\{X_{\overline 0}^{*}, Y_{\overline 0}^{*}, X_{\overline 1}^{*}, Y_{\overline 1}^{*}\}$ 
		      is the dual basic of 
		      $\{X_{\overline 0}, Y_{\overline 0}, X_{\overline 1}, Y_{\overline 1}\}$.
		\item
		 $H^2(\mathfrak{g}_{4,2}^s,\mathbb{C}) =\left\{ 0 \right\}$ .
		\item 
		$H^2(\mathfrak{g}_{6}^s,\mathbb{C}) = 
		  span\left\{ 
       		    \left[ Y_{\overline 0}^{*}\otimes X_{\overline 1}^{*} \right],
		        \left[ Y_{\overline 0}^{*}\otimes Y_{\overline 1}^{*} \right],
		        \left[ {\left( Z_{\overline 1}^{*} \right)}^2 \right],
		        \left[ {\left( T_{\overline 1}^{*} \right)}^2 \right],
		     \right.$
		\begin{flushright}
		$    \left. 
		        \left[ X_{\overline 1}^{*}Z_{\overline 1}^{*} 
		              - X_{\overline 0}^{*}\wedge Y_{\overline 0}^{*}\right],   
		        \left[ Y_{\overline 1}^{*}T_{\overline 1}^{*} 
		              - X_{\overline 0}^{*}\wedge Y_{\overline 0}^{*} \right] 
     		\right\}$  
      \end{flushright}		
		where 
		$\{X_{\overline 0}^{*}, Y_{\overline 0}^{*}, X_{\overline 1}^{*}, Y_{\overline 1}^{*}, Z_{\overline 1}^{*}, T_{\overline 1}^{*}\}$ is the dual basic of $\{X_{\overline 0}, Y_{\overline 0}, X_{\overline 1}, Y_{\overline 1}, Z_{\overline 1}, T_{\overline 1}\}$.
 \end{enumerate}
\end{thm}

\vskip0.3cm

\noindent {\bf The Proof of Theorem \ref{secondcohomology}}

\begin{enumerate} [(i)]
	\item 
	Firstly, we consider $\mathfrak{g}_{4,1}^s = \left( \mathbb{C}X_{\overline 0}\oplus \mathbb{C}Y_{\overline 0} \right)\oplus \left( \mathbb{C}X_{\overline 1}\oplus \mathbb{C}Y_{\overline 1} \right)$, where
	 $\mathfrak{g}_{\overline 0} = span \{X_{\overline 0}, Y_{\overline 0} \},$ $\mathfrak{g}_{\overline 1} = span \{X_{\overline 1}, Y_{\overline 1} \}.$
	 
	According to the paper \cite{DU14}, 
	the associated 3-form  of  $\mathfrak{g}_{4,1}^{s}$  is 
	$ I = Y_{\overline 0}^{*}\otimes {\left( Y_{\overline 1}^{*} \right)}^{2}$. 
	By a straightforward computation, we obtain
	
	\noindent $\left\{ I,X_{\overline 0}^{*} \right\}=  {\left( Y_{\overline 1}^{*} \right)}^{2}$, 
	          $\left\{ I,Y_{\overline 0}^{*} \right\} = 0$, 
	          $\left\{ I,X_{\overline 1}^{*} \right\}= 2Y_{\overline 0}^{*}\otimes Y_{\overline 1}^{*}$, 
	          $\left\{ I,Y_{\overline 1}^{*} \right\}=0$,
	
	\noindent $\left\{ I,X_{\overline 0}^{*}\wedge Y_{\overline 0}^{*} \right\}
	               = Y_{\overline 0}^{*}\otimes {\left( Y_{\overline 1}^{*} \right)}^{2}$, 
	          $\left\{ I,X_{\overline 0}^{*}\otimes X_{\overline 1}^{*} \right\}
	               =  X_{\overline 1}^{*} {\left( Y_{\overline 1}^{*} \right)}^{2} 
	                  + 2X_{\overline 0}^{*} \wedge Y_{\overline 0}^{*}\otimes Y_{\overline 1}^{*}$,
	
	\noindent $\left\{ I,X_{\overline 0}^{*} \otimes Y_{\overline 1}^{*} \right\}
	              = - {\left( Y_{\overline 1}^{*} \right)}^{3}$, 
	          $\left\{ I,Y_{\overline 0}^{*}\otimes X_{\overline 1}^{*} \right\}= 0$, 
	          $\left\{ I,Y_{\overline 0}^{*}\otimes Y_{\overline 1}^{*} \right\}= 0$, 
	
	\noindent $\left\{ I, {\left( X_{\overline 1}^{*} \right)}^{2} \right\}
	              =  4Y_{\overline 0}^{*}\otimes X_{\overline 1}^{*} Y_{\overline 1}^{*}$, 
	          $\left\{ I, {\left( Y_{\overline 1}^{*} \right)}^{2} \right\}=0$, 
	          $\left\{I, X_{\overline 1}^{*} Y_{\overline 1}^{*} \right\}
	              =  2Y_{\overline 0}^{*}\otimes {\left( Y_{\overline 1}^{*} \right)}^{2}$.
	
	Then we get 
	$$Im \delta_{1}=span\left\{ {\left( Y_{\overline 1}^{*} \right)}^{2}, Y_{\overline 0}^{*} \otimes Y_{\overline 1}^{*} \right\}$$ 
	and 
	$$Ker \delta _{2} = span\left\{ Y_{\overline 0}^{*} \otimes X_{\overline 1}^{*}, Y_{\overline 0}^{*}\otimes Y_{\overline 1}^{*}, {\left( Y_{\overline 1}^{*} \right)}^{2}, X_{\overline 1}^{*} Y_{\overline 1}^{*} - 2X_{\overline 0}^{*} \wedge Y_{\overline 0}^{*} \right\}.$$ 
	Therefore
	$$H^{2} (\mathfrak{g}_{4,1}^{s},\mathbb{C})=Ker \delta _{2}/ Im \delta _{1}
    	= span\left\{ 
    	      \left[ Y_{\overline 0}^{*} \otimes X_{\overline 1}^{*} \right], 
    	      \left[  X_{\overline 1}^{*}Y_{\overline 1}^{*} 
    	               - 2X_{\overline 0}^{*} \wedge Y_{\overline0}^{*} \right]
    	  \right\}.$$ 
\item
	 Next, we consider $\mathfrak{g}_{4,2}^{s}$ 
	$=\left( \mathbb{C} X_{\overline 0} \oplus \mathbb{C} Y_{\overline 0} \right)\oplus \left( \mathbb{C}X_{\overline 1} \oplus \mathbb{C} Y_{\overline 1} \right)$,
	where $ \mathfrak{g}_{\overline 0}$ 
	$=span\{ X_{\overline 0}, Y_{\overline 0} \}$ and 
	$ \mathfrak{g}_{\overline 1}$ 
	$=span\{ X_{\overline 1}, Y_{\overline 1} \}.$
	
	From \cite{DU14}, we obtain the associated 3-form $I = Y_{\overline 0}^{*} \otimes X_{\overline 1}^{*} Y_{\overline 1}^{*}$. By a similar computation as above, we have 
	$$ Ker \delta _{2} 
	  = Im \delta _{1} 
	  = span\left\{ X_{\overline 1}^{*} Y_{\overline 1}^{*}, 
	                Y_{\overline 0}^{*} \otimes X_{\overline 1}^{*}, 
	                Y_{\overline 0}^{*} \otimes Y_{\overline 1}^{*} 
	       \right\}.$$ 
Therefore we get 
$H^{2}(\mathfrak{g}_{4,2}^{s},\mathbb{C})=\left\{ 0 \right\}$.
	
	\item
	Finally, we consider {$\mathfrak{g}_{6}^{s}$ 
	$=\left( \mathbb{C} X_{\overline 0} \oplus \mathbb{C} Y_{\overline 0} \right)\oplus \left( \mathbb{C}X_{\overline 1} \oplus \mathbb{C} Y_{\overline 1} \oplus \mathbb{C} Z_{\overline 1} \oplus \mathbb{C} T_{\overline 1} \right)$,}
	where $\mathfrak{g}_{\overline 0} = span\{ X_{\overline 0}, Y_{\overline 0} \}   $, 
	$ \mathfrak{g}_{\overline 1} = span\{ X_{\overline 1}, Y_{\overline 1}, Z_{\overline 1}, T_{\overline 1}\}$.
	
	By a similar computation, we have 
	$$I = Y_{\overline 0}^{*} \otimes Z_{\overline 1}^{*} T_{\overline 1}^{*};
	Im \delta _{1} = span\left\{ Z_{\overline 1}^{*} T_{\overline 1}^{*}, Y_{\overline 0}^{*} \otimes T_{\overline 1}^{*}, Y_{\overline 0}^{*}\otimes Z_{\overline 1}^{*} \right\};$$
\begin{multline*}
	Ker \delta _{2} = span\left\{ Y_{\overline 0}^{*} \otimes X_{\overline 1}^{*}, Y_{\overline 0}^{*}\otimes Y_{\overline 1}^{*}, Y_{\overline 0}^{*} \otimes Z_{\overline 1}^{*},Y_{\overline 0}^{*}\otimes T_{\overline 1}^{*},{\left( Z_{\overline 1}^{*} \right)}^{2}, {\left( T_{\overline 1}^{*} \right)}^{2}, \right.\\ 
	\left. Z_{\overline{1}}^{*}T_{\overline{1}}^{*},X_{\overline{1}}^{*}Z_{\overline{1}}^{*}-X_{\overline{0}}^{*}\wedge Y_{\overline{0}}^{*},Y_{\overline{1}}^{*}T_{\overline{1}}^{*}-X_{\overline{0}}^{*}\wedge Y_{\overline{0}}^{*} \right\}.
		\end{multline*}
	
	Therefore we get
	\begin{multline*}
H^{2}(\mathfrak{g}_{6}^{s},\mathbb{C})=Ker \delta _2 / Im \delta _1		 
		  = span\left\{ 
       		    \left[ Y_{\overline 0}^{*}\otimes X_{\overline 1}^{*} \right],
		        \left[ Y_{\overline 0}^{*}\otimes Y_{\overline 1}^{*} \right],
		        \left[ {\left( Z_{\overline 1}^{*} \right)}^2 \right],
		     \right.\\
		     \left.
		        \left[ {\left( T_{\overline 1}^{*} \right)}^2 \right], 
		        \left[ X_{\overline 1}^{*}Z_{\overline 1}^{*} 
		              - X_{\overline 0}^{*}\wedge Y_{\overline 0}^{*}\right],   
		        \left[ Y_{\overline 1}^{*}T_{\overline 1}^{*} 
		              - X_{\overline 0}^{*}\wedge Y_{\overline 0}^{*} \right] 
     		\right\}.  
      \end{multline*}
\end{enumerate}
The proof is complete. $\hfill {\Box}$

\section{Classification of 8-dimensional Solvable \\Quadratic Lie Superalgebras Having \\6-dimensional Indecomposable Even Part}

The main result of the paper is a classification of all solvable quadratic Lie superalgebras of dimension 8 having indecomposable even part of dimension 6. 
In order to get this classification, it is necessary to observe consequence of adjoint orbits of the Lie algebra $\mathfrak{s}\mathfrak{p}(2)$ and the double extension.

\subsection{Adjoint Orbits of Symplectic Lie Algebra $\mathfrak{s}\mathfrak{p}(2)$}

Every non-zero element $\left( \begin{matrix}

   a & b  \\
   c & -a  \\
\end{matrix} \right)$ in the Lie algebra $\mathfrak{s}\mathfrak{p}(2)$ is either nilpotent or semisimple. It means that we can choose a basis of ${{\mathbb{C}}^{2}}$ such that the matrix representation of an element in $\mathfrak{s}\mathfrak{p}(2)$ with respect to this basis is 
diagonal or strictly upper triangliar.
  
\begin{lem} \label{lemma1} Let $A,B\in \mathfrak{s}\mathfrak{p}(2)$  such that $[A,B] = 0$. Then  $A$ and  $B$ are linearly dependent.
\end{lem}

\noindent {\it Proof.}
The Lie algebra $\mathfrak{sp}(2)$ has a basis $\left\{ H, X, Y \right\}$ with 
$$[H,X] = 2X,[H,Y] = -2Y,[X,Y] = H.$$ 
Set 
$$A = aH + bX + cY,B = a'H + b'X + c'Y.$$ 
From $\left[ {A,B} \right] = 0$, it implies 
$$(bc' - b'c)H + 2(ab' - a'b)X - 2(ac' - a'c)Y = 0.$$ That means $A$ and  $B$ are linearly dependent. \hfill $\square$

As an immediate consequence of Lemma \ref{lemma1}, we have the following lemma.

\begin{lem} \label{lemma2} 
 Let $A,B,C\in \mathfrak{s}\mathfrak{p}(2)$. If  $[A,B] = C$ and $[A,C] = [B,C] = 0$ then $C = 0$. \hfill $\square$ 
\end{lem}

\begin{lem} \label{lemma3}
 Let $A,B \in \mathfrak{s}\mathfrak{p}(2)$ such that $B \ne 0$ and $[A,B] = B$. Then $A$ is the semisimple element with eigenvalues  $\dfrac{1}{2}, - \dfrac{1}{2}$ and $B$ is nilpotent.
\end{lem}

\noindent {\it Proof.}
We can choose a basis such that $A$ is 
 $$\left( \begin{matrix}
   0 & 1  \\
   0 & 0  \\
\end{matrix} \right) \text{ or } 
\left( \begin{matrix}
   \lambda  & 0  \\
   0 & -\lambda   \\
\end{matrix} \right).$$
 If
 $A=\left( \begin{matrix}
   0 & 1  \\
   0 & 0  \\
\end{matrix} \right)=X$, let $B=aH+bX+cY$. Since $[A,B]=B$ one has 
$$-2aX+cH=aH+bX+cY.$$ 
So $B=0$ and this is a contradiction. Thus,  $A=\left( \begin{matrix}
   \lambda  & 0  \\
   0 & -\lambda   \\
\end{matrix} \right)=\lambda H$ . Set $B=aH+bX+cY$. By $[A,B]=B$, we have 
$$2\lambda bX-2c\lambda Y=aH+bX+cY.$$
This implies 
$$a=0, 2\lambda b=b, -2c\lambda =c.$$ Then  
$$a=0,\lambda =\dfrac{1}{2},c=0 \text{ or }a=0,\lambda =-\dfrac{1}{2},b=0.$$ 
The proof is complete. \hfill $\square$
\subsection{Double extension of quadratic Lie superalgebras}

In this subsection, we recall the notion of a double extension of quadratic Lie superalgebras which is introduced in \cite{BBB} and some its properties because it will be used later. For details we refer the reader to the paper \cite{Kac77} of V. G. Kac,  the paper \cite{BBB} of I. Bajo, S. Benayadi, M. Bordemann and the paper \cite{BB97} of I. Bajo, S. Benayadi. Firstly, we recall the definition of homogeneous superderivations.

\begin{defn}[see \cite{Kac77}]
	Let $\mathfrak{g}$ be a Lie superalgebra. An endomorphism $D \in Hom_{\alpha}(\mathfrak{g}, \mathfrak{g})$ (where $\alpha \in \mathbb{Z}_{2}$) is called {\it a  homogeneous superderivation} of degree $\alpha $ of $\mathfrak{g}$ if 
 $$D[X,Y]= [DX,Y]+(-1)^{\alpha x} [X,DY]; \, \forall x \in \mathbb{Z}_{2}, \, \forall X\in \mathfrak{g}_{x},\, \forall Y\in \mathfrak{g}.$$
\end{defn}
Denote by
$\left( Der(\mathfrak{g})\right)_{\alpha }\subset Hom_{\alpha}(\mathfrak{g}, \mathfrak{g})$ 
the space of all homogeneous superderivations of degree $\alpha $. Assuming that 
$Der(\mathfrak{g}) = \left( Der(\mathfrak{g}) \right)_{\overline 0} \oplus \left( Der(\mathfrak{g})\right)_{\overline 1}$, 
it is easily seen that $Der(\mathfrak{g})$ is Lie subsuperalgebra of Lie superalgebra  
$Hom(\mathfrak{g}, \mathfrak{g})$ 
and we call it the Lie superalgebra of
superderivations of $\mathfrak{g}$.
\begin{defn}[see \cite{BBB}]
Let $(\mathfrak{g},B)$ be a quadratic Lie superalgebra. A homogeneous superderivation $D$ of degree $\alpha $ of $\mathfrak{g}$ is called {\it skew-supersymmetric} if 
$$B(D(X),Y)=-{{(-1)}^{\alpha x}}B(X,D(Y)); \, \forall x \in \mathbb{Z}_{2}, \, \forall X\in \mathfrak{g}_{x},\, \forall Y\in \mathfrak{g}.$$
\end{defn}

It is proved that the vector subspace of $Der(\mathfrak{g})$ generated by the set of all
homogeneous skew-supersymmetric superderivations of $(\mathfrak{g}, B)$ is a Lie subsuperalgebra of $Der(\mathfrak{g})$ and it is denoted by $Der_{a}(\mathfrak{g},B)$.
In the remainder of this subsection we give the notion of a double extension of quadratic Lie superalgebras.

\begin{prop}[see \cite{BBB}, Theorem 2.4] \label{prop5}
 Let $\left( \mathfrak{g},B \right)$ be a quadratic Lie superalgebra and  $\mathfrak{h}$ be a Lie superalgebra. Suppose that $\psi :\mathfrak{h}\to Der_{a} \left( \mathfrak{g},B \right)$ is a morphism of Lie superalgebras. Define a bilinear map $\varphi$  from $\mathfrak{g}\times \mathfrak{g}$ to ${{\mathfrak{h}}^{*}}$ by 
$$\varphi (X,Y)(Z)=(-1)^{(x+y)z} B(\psi (Z)(X),Y); \, \forall x, y, z \in \mathbb{Z}_{2}, \forall X\in \mathfrak{g}_{x}, Y\in \mathfrak{g}_{y}, Z\in \mathfrak{h}_{z}.$$
Let $\pi _{\mathfrak{h}}$ by the coadjoint representation of $\mathfrak{h}$. Then the $\mathbb{Z}_{2}-$graded vector space $\overline{\mathfrak{g}} = \mathfrak{h} \oplus \mathfrak{g} \oplus \mathfrak{h}^{*}$ becomes a Lie superalgebra with the bracket defined by
\begin{multline*}
[Z+X+f,W+Y+g]= [Z,W]_{\mathfrak{h}} + [X,Y]_{\mathfrak{g}} + \psi (Z)(Y)-  (-1)^{xy}\psi (W)(X) \\
+\pi (Z)(g)-(-1)^{xy} \pi (W)(f)+\varphi (X,Y)
\end{multline*}
for all 
$Z+X+f\in \overline{\mathfrak{g}}_{x}, W+Y+g\in \overline{\mathfrak{g}}_{y}; \, x, y, z \in \mathbb{Z}_{2}$.
\noindent Furthermore, let $\gamma $ be an even supersymmetric invariant bilinear form on $\mathfrak{h}$. Then  $\overline{\mathfrak{g}}$ becomes a quadratic Lie superalgebra with the bilinear form defined by 
$$\overline{B}(Z+X+f,W+Y+g)=B(X,Y)+\gamma (Z,W)+f(W)+ (-1)^{xy}g(Z).$$  \hfill $\square$
\end{prop}
The quadratic Lie superalgebra \,$\overline{\mathfrak{g}}$ in Proposition 3.2.3 is called a double extension of  
$\left( \mathfrak{g}, B \right)$ by  $\mathfrak{h}$ (by means of $\psi $).

In particular, if 
$\mathfrak{h}= \mathbb{C}e$, $\mathfrak{h}^{*}  =\mathbb{C}f$ 
then $D=\psi (e)$ is an even skew-symmetric superderivation on $\mathfrak{g}$ and  
$\overline{\mathfrak{g}} = \mathbb{C}e\oplus \mathfrak{g}\oplus \mathbb{C}f$. In this case, $\overline{\mathfrak{g}}$ 
is called a {\it 1-dimensional double extension} of $\mathfrak{g}$ by means of $D$. The Lie super bracket on $\overline{\mathfrak{g}}$ is defined as follows
$$[X,Y]_{\overline{\mathfrak{g}}} = [X,Y]_{\mathfrak{g}} + (D(X),Y)f,[e,X]=D(X),[f,\overline{\mathfrak{g}}] = 0; \, \forall X,Y\in \mathfrak{g}.$$
The invariant bilinear form $\overline{B}$ on $\overline{\mathfrak{g}}$ is defined by
$$\overline{B}(e,f)=1,\overline{B}(e,\mathfrak{g})=\overline{B}(f,\mathfrak{g})=0,\overline{B}(X,Y)=B(X,Y),\forall X,Y\in \mathfrak{g}.$$

The following proposition presents the crucial role of one-dimensional double extensions for the construction of quadratic Lie superalgebras.

\begin{prop}[see \cite{BB97}] \label{prop6} 
Let $(\mathfrak{g},B)$  be an irreducible quadratic Lie superalgebra of dimension $n$ such that $n > 1$.  If $\mathfrak{z}(\mathfrak{g})\cap {{\mathfrak{g}}_{\overline{0}}}\ne \{0\}$ then $\mathfrak{g}$  is a double extension of a quadratic Lie superalgebra of dimension $n - 2$ by the one-dimensional Lie algebra.  \hfill$\square$ 
\end{prop}

\subsection{Classification of 8-dimensional Solvable Quadratic Lie Superalgebras having 6-dimensional Indecomposable Even Part}

The main result of this section is the classification of 8-dimensional solvable quadratic Lie algebras having 6-dimensional indecomposable even part. Before starting the main theorem, we first recall the classification of the 6-dimensional quadratic Lie algebras in \cite{PDL12}. 

\begin{prop}[see \cite{PDL12}] \label{prop7}
	 Let $( \mathfrak{g},B)$ be a solvable quadratic Lie algebra of dimension 6. Assume $ \mathfrak{g}$ indecomposable . Then there exists a basis $\{ Z_1, Z_2, Z_3, X_1, X_2, X_3 \}$ of $\mathfrak{g}$ such that the bilinear form $B$ is defined by   $B(X_i, Z_j)= \delta _{i,j}; 1 \le i,j \le 3$, the other are zero and $ \mathfrak{g}$ is i-isomorphic to each of Lie algebras as follows:
\begin{itemize} 
\item [(1)]

	  $\mathfrak{g}_{6,1}:$
	  $\left[X_1, X_2 \right] = Z_3, \left[X_2, X_3 \right] = Z_1, \left[X_3, X_1 \right] = Z_2.$
\item [(2)]
    $\mathfrak{g}_{6,2}(\lambda):$
	$\left[ X_3, Z_1 \right] = Z_1,
    \left[ X_3, Z_2 \right]=\lambda Z_2,
    \left[ X_3, X_1 \right] = -X_1,$
    $\left[ X_3, X_2 \right] = -\lambda X_2,
    \left[ Z_1, X_1 \right] = Z_3,
    \left[ Z_2, X_2 \right] = \lambda Z_3 $ 
    where 
    $\lambda \in \mathbb{C}$ and $\lambda \ne 0$.
    In this case 
    $\mathfrak{g}_{6,2}(\lambda _1)$ and $\mathfrak{g}_{6,2}(\lambda _2)$ 
    is i-isomorphic if and only if 
    $\lambda _2 = \pm \lambda _1$ or $\lambda _2 = \lambda _1^{-1}$. 
\item [(3)]
    $\mathfrak{g}_{6,3}:$
    $\left[ X_3, Z_1 \right]= Z_1, \left[ X_3, Z_2 \right]= Z_1 + Z_2, 
    \left[ X_3, X_1 \right]= -X_1 - X_2,$
    $\left[ X_3,X_2 \right]= - X_2,$
 
    $\left[ Z_1, X_1 \right]=\left[ Z_2, X_1 \right]=\left[ Z_2, X_2 \right]= Z_3.$ 
\end{itemize}
\end{prop}
\hfill $\square$

Now, we consider an 8-dimensional solvable quadratic Lie superalgebras 
$\left( \mathfrak{g} = \mathfrak {g}_{\overline 0}\oplus \mathfrak{g}_{\overline 1},B \right)$
having 6-dimensional indecomposable even part. Because $\mathfrak{g}$ is solvable, so is $\mathfrak{g}_{\overline 0}$. Since the even part $\mathfrak{g}_{\overline 0}$ is indecomposable, $\mathfrak{g}_{\overline 0}$ is isometrically isomorphic to each of quadratic Lie superalgebras $\mathfrak{g}_{6,1}, \mathfrak{g}_{6,2}(\lambda ), \mathfrak{g}_{6,3}$ which have been listed in Proposition \ref{prop7}.

Assume that $\mathfrak{g}$ is decomposable. Consider an arbitrary proper non-degenerate ideal $\mathfrak{j} = \mathfrak{j}_{\overline 0}\oplus \mathfrak{j}_{\overline 1}$ of $\mathfrak{g}$. Then $\mathfrak{j}^{\bot}$ is also a proper non-degenerate ideal of $\mathfrak{g}$. Hence,  $\mathfrak{j}_{\overline 0}$ or ${\left( \mathfrak{j}^{\bot} \right)}_{\overline 0}$ is a non-zero and non-degenerate ideal of $\mathfrak{g}_{\overline 0}$. Because $\mathfrak{g}_{\overline 0}$ is indecomposable, then  $\mathfrak{j}_{\overline 0} = \mathfrak{g}_{\overline 0}$ or ${\left( \mathfrak{j}^{\bot} \right)}_{\overline 0} = \mathfrak{g}_{\overline 0}$.  Therefore, without loss of generality, we may assume that $\mathfrak{j} = \mathfrak{g}_{\overline 0}\oplus \mathfrak{j}_{\overline 1}$. Then we get
\begin{center}
	$\mathfrak{j}_{\overline 1} = \{ 0 \} $, ${\left( {\mathfrak j}^{\bot} \right)}_{\overline 0}= \{ 0\}$ and ${\left({\mathfrak j}^{\bot} \right)}_{\overline 1} = {\mathfrak g}_{\overline 1}$.
\end{center}
That means that  
$\mathfrak{g} = {\mathfrak g}_{\overline 0} \mathop{\oplus}\, {\mathfrak g}_{\overline 1}$ 
where
${\mathfrak g}_{\overline 0}$ is an 6-dimensional indecomposable quadratic Lie algebra, ${\mathfrak g}_{\overline 1}$ is a symplectic vector space of dimension 2 and $[{\mathfrak g}_{\overline 0}, {\mathfrak g}_{\overline 1}]=\{ 0 \}$.

Next, we consider $\mathfrak{g}$  as an indecomposable quadratic Lie superalgebra with the even and odd parts are given as follows 
$${\mathfrak g}_{\overline 0} = span \{ X_i, Z_i \}, B\left( X_i, X_j \right) = B\left(Z_i, Z_j \right) = 0, B\left(X_i, Z_j \right) = \delta_{ij};\, 1\le i, j \le 3;$$
$${{\mathfrak{g}}_{\overline{1}}}=span\left\{ {{Y}_{\overline{1}}},{{T}_{\overline{1}}} \right\}, B\left( {{Y}_{\overline{1}}},{{T}_{\overline{1}}} \right)=1.$$

\vskip0.2cm
By Proposition \ref{prop7}, there are three cases for $\mathfrak{g}_{\overline 0}$.  Now, we will consider these cases one by one.

\begin{itemize}
	\item[(1)] {${{\mathfrak{g}}_{\overline{0}}}={{\mathfrak{g}}_{6,1}}$}
	
	In this case, the Lie brackets are defined by
	$$\left[X_1, X_2 \right] = Z_3, \, \left[X_2, X_3 \right] = Z_1, \, \left[X_3, X_1 \right] = Z_2. $$
	 
	Remark that 
	$$ad(X){\rvert}_{{\mathfrak g}_{\overline 1}} \in \mathfrak{s}\mathfrak{p}\left( {\mathfrak g}_{\overline 1}, B{\lvert}_{{\mathfrak g}_{\overline 1}\times {\mathfrak g}_{\overline 1}} \right)\cong \mathfrak{s}\mathfrak{p}\left( 2 \right), \, \forall X \in {\mathfrak g}_{\overline 0}.$$
	 
	By Lemma \ref{lemma2}, we get   
	$$ad(Z_1){\rvert}_{{\mathfrak g}_{\overline 1}}= ad(Z_2){\rvert} _{{\mathfrak g}_{\overline 1}} = ad(Z_3){\rvert}_{{\mathfrak g}_{\overline 1}} = 0.$$
	By Lemma \ref{lemma1}, any pair of elements in the set 
	$\{ad(X_1){\rvert}_{{\mathfrak g}_{\overline 1}}, \, ad(X_2){\rvert} _{{\mathfrak g}_{\overline 1}}, \, ad(X_3){\rvert}_{{\mathfrak g}_{\overline 1}}\}$ 
	is linearly dependent.
	Because $\mathfrak{g}$ is indecomposable, 
	one of $ad(X_1){\rvert}_{{\mathfrak g}_{\overline 1}}, \, ad(X_2){\rvert} _{{\mathfrak g}_{\overline 1}}$ and $ad(X_3){\rvert}_{{\mathfrak g}_{\overline 1}}$ 
	is non-zero. 
	On the other hand, by the classification of 
	$Sp\left({\mathfrak g}_{\overline 1} \right)$-orbits of the Lie algebra 
	$\mathfrak{s}\mathfrak{p}\left( {\mathfrak g}_{\overline 1} \right)$, 
	we can choose a Darboux basic $\left\{ Y_{\overline 1}, T_{\overline 1} \right\}$ 
	such that the representation matrices of  
	$ad(X_1){\rvert}_{{\mathfrak g}_{\overline 1}}, \,ad(X_2){\rvert} _{{\mathfrak g}_{\overline 1}}, \, ad(X_3){\rvert}_{{\mathfrak g}_{\overline 1}}$ 
	either are concurrently diagonal or strictly upper triangliar. Therefore, we need consider two different cases.
	
	\noindent 1a)
	$ad(X_1){\rvert}_{{\mathfrak g}_{\overline 1}} = \left( \begin{matrix}
	\lambda  & 0  \\
	0 & -\lambda   \\
	\end{matrix} \right),$ 
	$ad(X_2){\rvert}_{{\mathfrak g}_{\overline 1}} = \left( \begin{matrix}
	\mu  & 0  \\
	0 & -\mu   \\
	\end{matrix} \right),$  
	$ad(X_3){\rvert}_{{\mathfrak g}_{\overline 1}} = \left( \begin{matrix}
	\nu  & 0  \\
	0 & -\nu   \\
	\end{matrix} \right)$,
	 
	where $\lambda ,\mu , \nu \in \mathbb{C}$ such that 
	they are not simultaneously zero. In this case, we have
	 
	$$\left[X_1, Y_{\overline 1} \right] = \lambda Y_{\overline 1}, \, \left[ X_1, T_{\overline 1} \right] = -\lambda T_{\overline 1}, \, \left[X_2, Y_{\overline 1} \right] = \mu Y_{\overline 1},$$ 
	
	$$\left[X_2, T_{\overline 1} \right] = -\mu T_{\overline 1}, \, 
	\left[X_3, Y_{\overline 1} \right] = \nu Y_{\overline 1}, \,
	\left[X_3, T_{\overline 1} \right]= -\nu T_{\overline 1}.$$
	
	Because $B$ is invariant, non-degenerate and even, 	we obtain 
	$$[Y_{\overline 1}, Y_{\overline 1}] = 0,
	[T_{\overline 1}, T_{\overline 1}] = 0,
	[Y_{\overline 1}, T_{\overline 1}] = \lambda Z_1 + \mu Z_2 + \nu Z_3.$$
	 
	We denote this Lie superalgebra by $\mathfrak{g}_{8,2,1}^{s}(\lambda ,\mu ,\nu )$ 
	where $\lambda ,\mu , \nu \in \mathbb{C}$ such that 
	they are not simultaneously zero. 
	
	\noindent 1b)
	$ad(X_3){\rvert}_{{\mathfrak g}_{\overline 1}} = \left( \begin{matrix}
	0 & 1  \\
	0 & 0  \\
	\end{matrix} \right),$
	$ad(X_1){\rvert}_{{\mathfrak g}_{\overline 1}} = \left( \begin{matrix}
	0 & \lambda   \\
	0 & 0  \\
	\end{matrix} \right),$ 
	$ad(X_2){\rvert}_{{\mathfrak g}_{\overline 1}} = \left( \begin{matrix}
	0 & \mu   \\
	0 & 0  \\
	\end{matrix} \right)$,
	
	where $\lambda ,\mu , \in \mathbb{C}$. In this case, we have 
	$$\left[ X_3,T_{\overline{1}} \right]= Y_{\overline{1}}, \,
	\left[ X_1,T_{\overline{1}} \right]=\lambda Y_{\overline{1}}, \,
	\left[ X_2,T_{\overline{1}} \right]=\mu Y_{\overline{1}}.$$ 
	Because $B$ is invariant, non-degenerate and even, we get
	$[ Y_{\overline{1}},\mathfrak{g}_{\overline{1}} ]=0$.
	In a similar way of argument to the one used in the case (1a), we have  
	$$[ T_{\overline{1}},T_{\overline{1}}]=\lambda Z_1 + \mu Z_2 + Z_3,$$ 
	and the others on $\mathfrak{g}_{\overline 1}$ are zero.
	We denote this Lie superalgebra by $\mathfrak{g}_{8,2,2}^{s}(\lambda ,\mu )$
	where $\lambda ,\mu , \in \mathbb{C}$. 

	\item[(2)] $\mathfrak{g}_{\overline{0}}=\mathfrak{g}_{6,2}(\lambda)$
	
	The Lie bracket is given by 
$$\left[ X_3, Z_1 \right] = Z_1,
\left[ X_3, Z_2 \right]=\lambda Z_2,
\left[ X_3, X_1 \right] = -X_1,$$
$$\left[ X_3, X_2 \right] = -\lambda X_2,
\left[ Z_1, X_1 \right] = Z_3,
\left[ Z_2, X_2 \right] = \lambda Z_3 $$  
for $\lambda \ne 0$. By Lemma \ref{lemma2},  $ad\left( Z_3 \right)\left| _{\mathfrak{g}_{\overline 1}} \right.=0$. Let 
$$V=span\left\{ ad(X_i)\left| _{\mathfrak{g}_{\overline 1}} \right.,ad(Z_j)\left| _{\mathfrak{g}_{\overline 1}} \right.:i=1,2,3;j=1,2 \right\}.$$
 Then $ad(X_3)\left| _{\mathfrak{g}_{\overline 1}} \right.\ne 0$ since $\mathfrak{g}$ is indecomposable. By Lemma \ref{lemma1}, we get $\dim V=1,2.$

\noindent 2a) Assume that $\dim V=1$

We get  
$$ad(X_1)\left| _{\mathfrak{g}_{\overline 1}} \right.=0,
ad(X_2)\left| _{\mathfrak{g}_{\overline 1}} \right.=0,ad(Z_1)\left| _{\mathfrak{g}_{\overline 1}} =0 \right.,ad(Z_2)\left| _{\mathfrak{g}_{\overline 1}} \right.=0$$ and 
$$ [ \mathfrak{g}_{\overline 1}, \mathfrak{g}_{\overline 1}] \subset
 {[ \mathfrak{g}_{\overline 0},\mathfrak{g}_{\overline 0} ]}^{\bot } =
 Z( \mathfrak{g}_{\overline 0})=\mathbb{C}{Z_3}.$$ 
According to Proposition \ref{prop5}, $\mathfrak{g}$ 
is a one-dimensional double extension of   
$\mathfrak{q}_{\overline 0} \oplus \mathfrak{g}_{\overline 1}$ 
(where $\mathfrak{q}_{\overline 0} = span\left\{ Z_1, Z_2, X_1, X_2 \right\}$) 
by means of  $ad(X_3) \in \mathfrak{o}( \mathfrak{q}_{\overline 0})\oplus \mathfrak{sp}( \mathfrak{g}_{\overline 1})$.
By applying the classification of $Sp\left(\mathfrak{g}_{\overline 1} \right)-$orbits of 
$\mathfrak{sp}\left( \mathfrak{g}_{\overline 1} \right)$, 
we have the following cases of $ad(X_3)$.
\begin{enumerate}
	\item [(i)]
	$ad(X_3)=\left( \begin{matrix}
	1 & 0 & 0 & 0 & 0 & 0  \\
	0 & \lambda  & 0 & 0 & 0 & 0  \\
	0 & 0 & -1 & 0 & 0 & 0  \\
	0 & 0 & 0 & -\lambda  & 0 & 0  \\
	0 & 0 & 0 & 0 & 0 & 1  \\
	0 & 0 & 0 & 0 & 0 & 0  \\
	\end{matrix} \right)$, 
	we denote this Lie superalgebra by $\mathfrak{g}_{8,2,3}^{s}(\lambda )$ 
	where $\lambda \in \mathbb{C}$ and $\lambda \ne 0$ . 
	\item [(ii)]
	$ad(X_3)=\left( \begin{matrix}
	1 & 0 & 0 & 0 & 0 & 0  \\
	0 & \lambda  & 0 & 0 & 0 & 0  \\
	0 & 0 & -1 & 0 & 0 & 0  \\
	0 & 0 & 0 & -\lambda  & 0 & 0  \\
	0 & 0 & 0 & 0 & \mu  & 0  \\
	0 & 0 & 0 & 0 & 0 & -\mu   \\
	\end{matrix} \right)$,
	we denote this Lie superalgebra by $\mathfrak{g}_{8,2,4}^{s}(\lambda ,\mu )$ where $\lambda ,\mu , \in \mathbb{C}$ and $\lambda \ne 0,\mu \ne 0$. 
\end{enumerate}
\noindent 2b) Assume that $\dim V=2$ 

Without loss of generality, we assume that
$ad(Z_1)\left| _{\mathfrak{g}_{\overline 1}} \right.\ne 0$. 
We have 
$$ad(X_1)\left| _{\mathfrak{g}_{\overline 1}} \right.=xad(Z_1)\left| _{\mathfrak{g}_{\overline 1}} \right.,$$
for some $x \in \mathbb{C}$ since  
$ad( X_1) \left| _{\mathfrak{g}_{\overline 1}} \right.,ad(Z_1)\left| _{\mathfrak{g}_{\overline 1}} \right.$ 
are linearly dependent. 

Moreover,
$$-ad(X_1) \left| _{\mathfrak{g}_{\overline 1}} \right.=\left[ ad(X_3) \left| _{\mathfrak{g}_{\overline 1}} \right., ad(X_1) \left| _{\mathfrak{g}_{\overline 1}} \right. \right]= xad(Z_1)\left| _{\mathfrak{g}_{\overline 1}} \right.$$ 
clearly force $ad(X_1) \left| _{\mathfrak{g}_{\overline 1}} \right.=0$.
$$\left[ ad(X_3) \left| _{\mathfrak{g}_{\overline 1}} \right., ad(Z_1) \left| _{\mathfrak{g}_{\overline 1}} \right. \right] = ad(Z_1) \left| _{\mathfrak{g}_{\overline 1}} \right.$$ 
infers that
$ad(X_3) \left| _{\mathfrak{g}_{\overline 1}} \right.$ 
is semisimple with eigenvalues $\dfrac{1}{2},-\dfrac{1}{2}$  and 
$ad(Z_1) \left| _{\mathfrak{g}_{\overline 1}} \right.$ 
is nilpotent (see Lemma \ref{lemma3}). Consequently, we can choose a Darboux basis
$\{ Y_{\overline 1}, T_{\overline 1} \}$ of  
$\mathfrak{g}_{\overline 1}$ such that 
$$ad(X_3)\left| _{\mathfrak{g}_{\overline 1}} \right.=\left( \begin{matrix}
\tfrac{1}{2} & 0  \\
0 & -\tfrac{1}{2}  \\
\end{matrix} \right),
ad(Z_{1})\left| _{\mathfrak{g}_{\overline 1}} \right.=\left( \begin{matrix}
0 & 1  \\
0 & 0  \\
\end{matrix} \right).$$
Assume that 
$$ad(X_2)\left| _{\mathfrak{g}_{\overline 1}} \right.=\alpha ad(Z_1)\left| _{\mathfrak{g}_{\overline 1}} \right., ad(Z_2)\left| _{\mathfrak{g}_{\overline 1}} \right.=\beta ad(Z_1)\left| _{\mathfrak{g}_{\overline 1}} \right.$$ 
for some $\alpha , \beta \in \mathbb{C}$. Then  
$$\left[ ad(X_3)\left| _{\mathfrak{g}_{\overline 1}}  \right., ad(Z_2)\left| _{\mathfrak{g}_{\overline 1}}  \right. \right] =\lambda ad(Z_2)\left| _{\mathfrak{g}_{\overline 1}}  \right.$$
implies that 
$\lambda =1$ or 
$ad(Z_2)\left| _{\mathfrak{g}_{\overline 1}} \right.=0$.
$$\left[ ad(X_3)\left| _{\mathfrak{g}_{\overline 1}} \right., ad(X_2)\left| _{\mathfrak{g}_{\overline 1}} \right. \right] =-\lambda ad(X_2)\left| _{\mathfrak{g}_{\overline 1}} \right.$$ 
yields that
$\lambda =-1$ or
$ad(X_2)\left| _{\mathfrak{g}_{\overline 1}} \right.=0$.

Note that when $\lambda =-1$, it can be returned to the case $\lambda =1$ by replacing  $X_2$ by $Z_2$ and $Z_2$ by $X_2$.

Therefore, we have two following cases.
\begin{enumerate}
	\item [(i)]
	If  $ad(Z_2)\left| _{\mathfrak{g}_{\overline 1}} \right.=ad(X_2)\left| _{\mathfrak{g}_{\overline 1}} \right.=0$, 
	then 
	$$[Y_{\overline 1},T_{\overline 1}] =\dfrac{1}{2} Z_3,[T_{\overline 1},T_{\overline 1}] = X_1.$$
	Denote this Lie superalgebra  by 
	$\mathfrak{g}_{8,2,5}^{s}(\lambda )$ where $\lambda \in \mathbb{C}$ and $\lambda \ne 0$.
	\item [(ii)]
	If  $\lambda =1$, then 
	$$ad(X_2)\left| _{\mathfrak{g}_{\overline 1}}  \right.=0, ad(Z_2)\left| _{\mathfrak{g}_{\overline 1}}  \right.= \left( \begin{matrix}
	0 & \mu   \\
	0 & 0  \\
	\end{matrix} \right),\mu \ne 0.$$
	By the invariance of the bilinear form 
	$B$, we obtain
	$$[Y_{\overline 1},T_{\overline 1}] =\dfrac{1}{2} Z_3,[T_{\overline 1}, T_{\overline 1}] =  X_1 +\mu X_2.$$ 
	Denote this Lie superalgebra by 
	$\mathfrak{g}_{8,2,6}^{s}(\lambda , \mu )$ 
	where $\lambda ,\mu , \in \mathbb{C}$ and $\lambda \ne 0,\mu \ne 0$. . 
\end{enumerate}	
		
	\item[(3)] $\mathfrak{g}_{\overline 0}= \mathfrak{g}_{6,3}$

The non-zero Lie brackets: 
$$\left[ X_3, Z_1 \right]= Z_1, \left[ X_3, Z_2 \right]= Z_1 + Z_2, \left[ X_3, X_1 \right]= -X_1 - X_2,$$
$$\left[ X_3,X_2 \right]= - X_2,
\left[ Z_1, X_1 \right]=\left[ Z_2, X_1 \right]=\left[ Z_2, X_2 \right]= Z_3.$$ 
Denote by
$$V=span\left\{ ad(X_1)\left| _{\mathfrak{g}_{\overline 1}} \right., ad(X_2)\left| _{\mathfrak{g}_{\overline 1}} \right., ad(X_3)\left| _{\mathfrak{g}_{\overline 1}} \right., ad(Z_1)\left| _{\mathfrak{g}_{\overline 1}} \right., ad(Z_2)\left| _{\mathfrak{g}_{\overline 1}} \right. \right\}.$$ 
Similarly to Case  (2), we get $ad(Z_3)\left| _{\mathfrak{g}_{\overline 1}} = 0 \right.$ and $\dim V=1,2.$

\noindent 3a) Assume that $\dim V=1$.

By the argument analogous to that used in Case (2), $\mathfrak{g}$  is a one-dimensional double extension of   
$\mathfrak{q}_{\overline 0} \oplus \mathfrak{g}_{\overline 1}$ 
(where $\mathfrak{q}_{\overline 0} = span\left\{ Z_1, Z_2,  X_1, X_2 \right\}$) by means of  $ad(X_3) \in \mathfrak{o}( \mathfrak{q}_{\overline 0}) \oplus \mathfrak{sp}( \mathfrak{g}_{\overline 1})$. Applying the classification of $Sp\left( \mathfrak{g}_{\overline 1} \right)-$orbits of $\mathfrak{sp}\left( \mathfrak{g}_{\overline  1} \right)$, we have the following cases of $ad(X_3)$. 
\begin{enumerate}
	\item [(i)]
	$ad(X_3)=\left( \begin{matrix}
	1 & 1 & 0 & 0 & 0 & 0  \\
	0 & 1 & 0 & 0 & 0 & 0  \\
	0 & 0 & -1 & 0 & 0 & 0  \\
	0 & 0 & -1 & -1 & 0 & 0  \\
	0 & 0 & 0 & 0 & 0 & 1  \\
	0 & 0 & 0 & 0 & 0 & 0  \\
	\end{matrix} \right)$,
	denote this Lie superalgebra by $\mathfrak{g}_{8,2,7}^{s}$.
	\item [(ii)]
	$ad(X_3)=\left( \begin{matrix}
	1 & 1 & 0 & 0 & 0 & 0  \\
	0 & 1 & 0 & 0 & 0 & 0  \\
	0 & 0 & -1 & 0 & 0 & 0  \\
	0 & 0 & -1 & -1 & 0 & 0  \\
	0 & 0 & 0 & 0 & \lambda  & 0  \\
	0 & 0 & 0 & 0 & 0 & -\lambda   \\
	\end{matrix} \right)$,
	denote this Lie superalgebra by $\mathfrak{g}_{8,2,8}^{s}(\lambda )$
	where 
	$0 \ne \lambda \in \mathbb{C}$.
\end{enumerate}
\noindent 3b) Assume that  $\dim V=2$. 

Firstly, we see that 
$ad(Z_2)\left| _{\mathfrak{g}_{\overline 1}}\right.$ 
or
$ad(X_1)\left| _{\mathfrak{g}_{\overline 1}} \right.$ 
is non-zero. 
There is no loss of generality in assuming 
$ad(Z_2)\left| _{\mathfrak{g}_{\overline 1}} \ne 0 \right.$. 
By Lemma \ref{lemma3}, we obtain 
$$ad(X_3)\left| _{\mathfrak{g}_{\overline 1}}  \right.=\left( \begin{matrix}
\tfrac{1}{2} & 0  \\
0 & -\tfrac{1}{2}  \\
\end{matrix} \right),
ad(Z_2)\left| _{\mathfrak{g}_{\overline 1}}  \right.=\left( \begin{matrix}
0 & 1  \\
0 & 0  \\
\end{matrix} \right),$$ 
$$ad(Z_1)\left| _{\mathfrak{g}_{\overline 1}}  \right.=
ad(X_1)\left| _{\mathfrak{g}_{\overline 1}}  \right.=
ad(X_2)\left| _{\mathfrak{g}_{\overline 1}}  \right.=0.$$
By the invariance of $B$, we have 
$$ [ X_3, Y_{\overline 1} ] =\dfrac{1}{2} Y_{\overline 1}, [ X_3, T_{\overline 1}] = - \dfrac{1}{2} T_{\overline 1},
 [Z_2, T_{\overline 1}] = Y_{\overline 1},
 [Y_{\overline 1}, T_{\overline 1}] = \dfrac{1}{2} Z_3,
  [T_{\overline 1}, T_{\overline 1}] = X_2.$$ 
Denote this Lie superalgera by $\mathfrak{g}_{8,2,9}^{s}$.	
\end{itemize}
Now, combining the above calculations, we obtain the main theorem of the paper  which gives the classification of 8-dimensional solvable quadratic Lie superalgebras having 6-dimensional indecomposable even part. 

\begin{thm}  Let $\mathfrak{g}$ be a solvable quadratic Lie superalgebra of dimension 8 such that its even part is indecomposable and 6-dimensional.
\begin{enumerate}
\item[(i)]
If  $\mathfrak{g}$ is decomposable then  $\mathfrak{g}$ is isomorphic to 
$\mathfrak{g}= \mathfrak{g}_{\overline 0} \overset{\bot }{\mathop{\oplus }}  \mathfrak{g}_{\overline 1}$ 
where 
$\mathfrak{g}_{\overline 0}$ 
is isometrically isomorphic to 
$\mathfrak{g}_{6,1},\mathfrak{g}_{6,2}(\lambda ),\mathfrak{g}_{6,3}$ 
which are given in Proposition 3.3.1; 
$\mathfrak{g}_{\overline 1}$ 
is a symplectic vector space of dimension 2 and
  $ [  \mathfrak{g}_{\overline 0},\mathfrak{g}_{\overline 1}] = \{ 0 \}$.
 \item[(ii)]
	If $\mathfrak{g}$ is indecomposable then $\mathfrak{g}$ is i-isomorphic to each of the following quadratic Lie superalgebras $\mathfrak{g}_{8,2,1}^{s}(\lambda ,\mu ,\nu )$, $\mathfrak{g}_{8,2,2}^{s}(\lambda ,\mu )$, $\mathfrak{g}_{8,2,3}^{s}(\lambda )$,
$\mathfrak{g}_{8,2,4}^{s}(\lambda ,\mu )$, $\mathfrak{g}_{8,2,5}^{s}(\lambda )$, $\mathfrak{g}_{8,2,6}^{s}(\lambda, \mu )$, $\mathfrak{g}_{8,2,7}^{s}$ , $\mathfrak{g}_{8,2,8}^{s}(\lambda )$ and $\mathfrak{g}_{8,2,9}^{s}$.
\hfill $\square$ 
\end{enumerate}
\end{thm}



\begin{thebibliography}{99}

\bibitem{BL17}
W. Bai, W. Liu, Cohomology of Heisenberg Lie Superalgebras, {\it J. of Math. Physics} {\bf 58}   (2017) 021701.
\bibitem{Bor97}
M. Bordemann, Nondegenerate invariant bilinear forms on nonassociative algebras, {\it Acta. Math. Uni. Comenianac} {\bf LXVI} (2) (1997) 151-201.
\bibitem{BB97}
I. Bajo, S. Benayadi, Lie algebras admitting a unique quadratic structure, {\it Commun. in Algebra} {\bf25} (9) (1997) 2795–2805.
\bibitem{BBB}
I. Bajo, S. Benayadi, M. Bordemann, Generalized double extension and descriptions of quadratic Lie superalgebras, arXiv:0712.0228v1 (2007). 
\bibitem{CD15}
T. T. H. Cao, M. T. Duong, The Betti numbers and the vector space of skew-symmetric derivations of solvable quadratic Lie algebras with dimension $\le $ 7, {\it  J. of Science}, Ho Chi Minh city University of Education, No {\bf 5} (70) (2015) 86-96 (in Vietnamese). 
\bibitem{Duo13}
M. T. Duong, The cohomology group ${{H}^{2}}(\mathfrak{g},\mathbb{C})$ of the elementary quadratic Lie algebras, {\it J. of Science}, Ho Chi Minh city University of Education, No {\bf 47} (81) (2013) 25 – 36.
\bibitem{DU14}
M. T. Duong, R. Ushirobira, Singular quadratic Lie superalgebras, {\it J. of Algebra} {\bf 407} (2014) 372–412.
\bibitem{FL84}
D.B.Fuchs, D.A.Leites,  Cohomology of Lie superalgebras, {\it Dokl.Bolg. Akad. Nauk} {\bf 37} (10) (1984) 1294-1296. 
\bibitem{Kac77}
V. G. Kac,    A sketch of Lie superalgebra theory, {\it Commun. Math. Physics}, {\bf 53}  (1977)  31—64.
\bibitem{MR85}
A. Medina, P. Revoy, Algèbres de Lie et produit scalaire invariant, {\it Ann. Sci. Éc. Norm. Sup., 4ème sér.} {\bf 18} (1985) 553-561.
\bibitem{MPU09}
I.A. Musson, G. Pinczon, R. Ushirobira, Hochschild cohomology and deformations of Clifford–Weyl algebras,   {\it SIGMA} {\bf 5} 27 (2009).
\bibitem{PDL12}
T. D. Pham, M. T. Duong, A. V. Le, Solvable quadratic Lie algebras in low dimensions, {\it East-West J. of Math.} {\bf 14} (2) (2012) 208-218.
\bibitem{PU07}
G. Pinczon, R. Ushirobira, New Applications of Graded Lie Algebras to Lie Algebras, Generalized Lie Algebras, and Cohomology, {\it J. Lie Theory}, {\bf 17} (2007) 633-667.
\bibitem{SZ98}
M.Scheunert, R.B.Zhang, Cohomology of Lie superalgebras and their generalizations, {\it J.Math. Phy.} {\bf 39} (9) (1998) 5024-5061.
\bibitem{SZ07}
Y.C.Su, R.B.Zhang, Cohomology of Lie superalgebras   $\mathfrak{s}{{\mathfrak{l}}_{m|n}}$ and $\mathfrak{o}\mathfrak{s}{{\mathfrak{p}}_{2|2n}}$, {\it Proc. London Math. Soc.} {\bf 94} (3) (2007) 91-136.

\end{thebibliography}
\end{document}